\documentclass[11pt,english,reqno]{amsart}
\usepackage[T1]{fontenc}
\usepackage{mathrsfs}
\usepackage{amsmath,amsthm,amssymb}
\usepackage{pxfonts}
\usepackage{latexsym}
\usepackage{times}
\usepackage{enumerate}
\usepackage{mathrsfs}
\usepackage{stmaryrd}
\usepackage{amsopn}
\usepackage{amsmath}
\usepackage{amssymb}
\usepackage{amsfonts}
\usepackage{amsbsy}
\usepackage{amscd,indentfirst,epsfig}
\usepackage{hyperref}
\usepackage{amsfonts,amsmath,latexsym,amssymb,verbatim,amsbsy}
\usepackage{amsthm}
\usepackage{dsfont}
\usepackage{colordvi}
\usepackage{pstricks}
\newtheorem{theo}{Theorem}[section]
\newtheorem{lemme}[theo]{Lemma}
\newtheorem{propo}[theo]{Proposition}
\newtheorem{cor}[theo]{Corollary}

\newtheorem{nb}[theo]{Remark}

\theoremstyle{definition}

\def \leq {\leqslant}
\def \geq {\geqslant}
\numberwithin{equation}{section}
\def\ind#1{\lower5pt\hbox{$\scriptstyle #1$}}

\def \x {x}
\def \d {{d}}
\def \e {\epsilon}
\def \e2 {\e^2}
\def \vt {(V(t))_{t \geq 0}}
\def \uot {(U(t))_{t \geq 0}}
\def \zt {(Z(t))_{t \geq 0}}

\def \esup {\operatornamewithlimits{ess\,sup}}

\def \l {\lambda}
\def \a {\alpha}
\def \b {\beta}
\def \v {v_{\star}}
\def \w {w_{\star}}
\def \D {\mathscr{D}}
\def \g {\overline{\gamma}}

\def\Q {{Q}}

\def\R{{\mathbb R}^3}
\def \S {{\mathbb S}^2}
\def \q {q}
\def \n {n}
\def \qn {|\q \cdot \n|}
\def \u {{u}_1}
\def \v {{v}}

\def \vb {\v_{\star}}
\def \w {{w}}

\def \wb {\w_{\star}}

\def \It {\int_{\R \times \S}}

\def \IRR {\int_{\R \times \R}}

\def\ep{\epsilon}

\title[Integral representation of the linear Boltzmann operator for granular gas]
{Integral representation of the linear Boltzmann operator for
granular gas dynamics with applications}

\author{L. Arlotti \& B. Lods}

\address{Luisa  Arlotti, Dipartimento di Ingegneria Civile, Universit\`a di
Udine, via delle Scienze 208,\newline \indent  33100 Udine,
Italy.}\email{luisa.arlotti@uniud.it}

\address{Bertrand Lods, Laboratoire de Math\'{e}matiques, cnrs umr 6620, Universit\'{e} Blaise Pascal\newline \indent (Clermont--Ferrand
2),
 63177 Aubi\`{e}re Cedex, France.}
\email{lods@math.univ-bpclermont.fr}
\begin{document}
\bibliographystyle{plain}
\maketitle
\begin{abstract} We investigate the properties of the collision operator $\Q$ associated to the linear Boltzmann equation
for {\it dissipative} hard-spheres arising in granular gas dynamics.
We establish that, as in the case of non--dissipative interactions,
the gain collision operator is an integral operator whose kernel is
made explicit. One deduces from this result a complete picture of
the spectrum of $\Q$ in an Hilbert space setting, generalizing
results from T. Carleman \cite{carleman} to granular gases. In the
same way, we obtain from this integral representation of $\Q^+$ that
the semigroup in $L^1(\R \times \R,\d \x \otimes \d\v)$ associated
to the linear Boltzmann equation for dissipative hard spheres is
{\it honest} generalizing known results from \cite{Ar}.

\end{abstract}
%

%
%
%
\section{Introduction}

We deal in this paper with the linear Boltzmann equation for
\textit{dissipative interactions} modeling the evolution of a
granular gas, undergoing {\it inelastic collisions} with its
 underlying medium. Actually, we shall see in the sequel that there
is no contrast  between the scattering theory of granular gases and
that of classical {\it (elastic)} gases. This may seem quite
surprising if one has in mind the fundamental differences that may
be emphasized between the {\it nonlinear kinetic theory} of granular
gases and that of classical gases, as briefly recalled in the next
lines.

\subsection{Granular gas dynamics: linear and nonlinear models}  Let us begin by recalling
the general features of the kinetic description of granular gas
dynamics that can be recovered from the monograph \cite{poschel} or
the more mathematically oriented survey \cite{villani}. If
$f(\x,\v,t)$ denotes the distribution function  of granular
particles with position $\x \in \R$ and velocity $\v \in \R$ at time
$t \geq 0,$ then the evolution of $f(\x,\v,t)$ is governed by the
following generalization of Boltzmann equation
\begin{equation}\label{BoltzNL}
\partial_t f(\x,\v,t)+\v \cdot \nabla_\x
f(\x,\v,t) = \mathcal{C}(f)(\x,\v,t),\end{equation} with initial
condition $f(\x,\v,0)=f_0(\x,\v) \in L^1(\R \times \R,\d x \otimes
\d \v)$, where the right--hand side $\mathcal{C}(f)$ models the
collision phenomena and depends on the phenomenon we describe.

In the {\bf nonlinear description}, the collision operator
$\mathcal{C}(f)=:\mathcal{B}[f,f]$ is a quadratic operator modeling
the binary collision phenomena between self-interacting particles.
For \textit{hard--spheres} interactions, it reads
\begin{equation*}\label{modelnonli}\mathcal{B}[f,f](\v)= \It \qn
\left[\dfrac{1}{\epsilon^2}f(\x,\vb,t)f(\x,\wb,t)-f(\x,\v,t)f(\x,\w,t)\right]\d\w
\d\n,\end{equation*} where  $\q$ is the relative velocity,
$\q=\v-\w.$ The microscopic velocities $(\vb,\wb)$ are the
pre--collisional velocities of the so--called inverse collision,
which results in $(\v,\w)$ as post-collisional velocities. The main
peculiarity of the kinetic description of granular gas is the
inelastic character of the microscopic collision mechanism  which
induces that \textit{part of the total kinetic energy is
dissipated}. This energy dissipation might be due to the roughness
of the surface or just to a non-perfect restitution and is measured
through a restitution coefficient $0 < \ep<1$ (which we assume here
to be constant for simplicity, see Remark \ref{referee3}). As a
consequence, the collision phenomenon is a \textit{non
microreversible process}. Generally, we assume that the energy
dissipation does not affect the conservation of momentum. Therefore,
in the homogeneous setting, i.e. when $f_0(\x,\v)=f_0(\v)$ is
independent of the position,  the number density of the gas  is
constant while the bulk velocity
  is conserved. However, the temperature of the gas
$$\vartheta(t)=\dfrac{1}{3}\int_{\R} |\v|^2 f(t,\v)\d\v$$
continuously decreases \textit{(cooling of granular gas)}. As a
consequence, the stationary state of the inelastic collision
operator $\mathcal{B}$ is a given Dirac mass. However, the
homogeneous Boltzmann equation for granular gases exhibits
self--similar solution (\textit{homogeneous cooling state})
\cite{ernst,mischler}. Note the important contrast with the
classical kinetic theory, i.e. when $\ep=1$, for which it is
well--known that the steady state of the collision operator is a
Maxwellian distribution.

The  {\bf linear Boltzmann equation for dissipative interactions}
concerns dilute particles (test particles with negligible mutual
interactions) immersed in a fluid at thermal equilibrium
\cite{loto,piasecki,spiga}. The total kinetic energy is dissipated
when the dilute particles collide with particles of the host fluid.
Such physical  models are well-suited to the study of the dynamics
of a mixture of impurities in a gas \cite{ferrari, pareschi} for
which the background is in thermodynamic equilibrium and that the
polluting particles are sufficiently few.  We refer the reader to
\cite{garzo} and the survey \cite{garzosurvey} for more details on
the theory of granular gaseous mixtures. Assuming the fluid at
thermal equilibrium and neglecting the mutual interactions of both
the test and dilute particles, the collision operator
$\mathcal{C}(f)=\Q(f)$ is a \textit{linear scattering operator}
given by
\begin{equation}\label{linearoperator}
\Q(f)=\mathcal{B}[f,\mathcal{M}_1]=\It \qn
\left[\dfrac{1}{\epsilon^2}f(\x,\vb,t)\mathcal{M}_1(\wb)-f(\x,\v,t)\mathcal{M}_1(\w)\right]\d\w
\d\n\end{equation} where $\mathcal{M}_1$ stands for the distribution
function of the host fluid. Note that in such a scattering model,
the microscopic masses of the dilute particles $m$ and that of the
host particles $m_1$ can be different. We will assume throughout
this paper that the distribution function of the host fluid is a
given normalized Maxwellian function:
$$\mathcal{M}_1(\v)=\bigg(\dfrac{m_1}{2\pi
\vartheta_1}\bigg)^{3/2}\exp
\left\{-\dfrac{m_1(\v-u_1)^2}{2\vartheta_1}\right\}, \qquad
\qquad \v \in \R,$$
where $u_1 \in \R$ is the given bulk velocity  and $\vartheta_1  >0$
is the given effective temperature of the host fluid. It can be
shown in this case that the number density  of the dilute gas is the
unique conserved macroscopic quantity (as in the elastic case). The
temperature is still not conserved but it remains bounded away from
zero, which prevents the solution to the linear Boltzmann equation
to converge towards a Dirac mass. This strongly contrasts to the
nonlinear description and suggests that the linear scattering model
associated with granular gases does not contrast too much with the
one associated with classical gases.

The first mathematical result in this direction is the following one
according to which, as in the classical case, the unique steady
state of $\Q$ remains Gaussian. The fact that the linear Boltzmann
equation still possesses a stationary Maxwellian velocity
distribution was first obtained in \cite{piasecki} and we refer to
\cite{loto} for a complete proof (existence and uniqueness) for
hard-spheres model (see also \cite{spiga} for a version of this
result for Maxwell molecules) :
\begin{theo}\label{existence} The  Maxwellian velocity distribution:
$$\mathcal{M}(\v)=\left(\dfrac{m}{2\pi
\vartheta^{\#}}\right)^{3/2}\exp
\left\{-\dfrac{m(\v-\u)^2}{2\vartheta^{\#}}\right\} \qquad \v \in
\R,$$ 
with $ \vartheta^{\#}=\dfrac{(1+\epsilon)
m}{2m+(3+\epsilon)m_1}\vartheta_1$ is the unique equilibrium state
of $\Q$ with unit mass.
\end{theo}
\begin{nb}\label{referee3}
Note that, if one does not assume the restitution coefficient
$\epsilon$ to be constant (see  \cite{poschel} for the general
expression of non-constant restitution coefficient
$\epsilon=\epsilon(\q)$ in the case, e.g., of visco-elastic spheres)
then the nature of the equilibrium state of $\Q$ is still an open
question: it is not known whether such a steady   state is a
Maxwellian or not. Consequently, it is still not clear that linear
inelastic scattering models behave like elastic ones. For this
reason, we shall restrict here our study of the linear Boltzmann
equation to a \textbf{\textit{constant}} restitution coefficient. We
also point out that, if the distribution function of the host fluid
$\mathcal{M}_1$ is not of gaussian type, the explicit expression of
the equilibrium state of $\Q$ is an open question to our
knowledge.\end{nb} The existence and uniqueness of such an
equilibrium state allows to establish a linear version of the famous
$H$--Theorem. Precisely, for any {\it convex} $C^1$--function $\Phi
\::\:\mathbb{R}^+ \to \mathbb{R}$, one can define the associated
entropy functional as
\begin{equation} \label{Hphi} H_{\Phi}(f|\mathcal{M})=\int_{\R}\mathcal{M}(\v)\,
\Phi\left (\dfrac{f(\v)}{\mathcal{M}(\v)} \right)\d\v , \end{equation}
\begin{theo}[$H$--Theorem \cite{loto,rolf}]\label{theoH} Let   $f_0(\v)$ be a space homogeneous distribution function with unit
mass and finite entropy, i.e.  $H_{\Phi}(f_0|\mathcal{M}) < \infty.$
Then,
\begin{equation}\label{h}
\dfrac{\d}{\d t}H_{\Phi}(f(t)|\mathcal{M})\leq 0 \qquad \qquad (t
\geq 0),
\end{equation}
where $f(t)$ stands for the (unique) solution to \eqref{BoltzNL} in
$L^1(\R,\d\v)$.\end{theo}
 Note that
such a result is valid for any scattering operator with positive
kernel and positive equilibrium \cite{rolf}. As an important
consequence, it can be shown by suitable compactness arguments that
any solution to the Boltzmann equation \eqref{BoltzNL} (with unit
mass) converges towards the Maxwellian equilibrium $\mathcal{M}$.
Note also that, for the nonlinear Boltzmann equation for dissipative
interactions, the temperature is a trivial Lyapunov functional
leading to the convergence of any solution towards a delta mass.
However, the construction of a Lyapunov functional in the
self-similar variables allowing relaxation towards the homogeneous
cooling state is still an open question (see, e.g. \cite{mischler}
for related problems).

To summarize, the steady state of the linear collision operator for
dissipative interactions has the same nature (a Maxwellian
distribution) as the one corresponding to non--dissipative
interactions. Moreover, as in the classical case, by virtue of the
$H$--Theorem, such a steady state attracts any solution to the space
homogeneous Boltzmann equation \eqref{BoltzNL}. This seems to
indicate that  \textit{most of the properties of the linear
Boltzmann equation for elastic interactions remain valid for
inelastic scattering models.}  It is the main subject of this paper
to make precise and confirm such an indication and the key
ingredient will be the derivation of an integral representation of
the gain part of the collision operator.

\subsection{Main results} The main concern of our paper is the derivation of a
suitable representation of the gain part of the collision operator
$\Q$ as an integral operator with\textit{ explicit kernel}.
Precisely, the linear collision operator $\Q$ can be split into
$\Q(f)=\Q^+(f)-\Q^-(f),$ where the gain part is
$$\Q^+(f)(\v)=\epsilon^{-2}\It |\q \cdot \n|f(\v_{\star})\mathcal{M}_1(\w_{\star})\d \w \d
\n$$ while
$$\Q^-(f)(\v)=\It |\q \cdot \n|f(\v)\mathcal{M}_1(\w)\d \w \d
\n= \sigma(\v)f(\v)$$ where the collision frequency $\sigma(\v)$ is
given by $\sigma(\v)=\It |\q \cdot \n|\mathcal{M}_1(\w)\d \w\d \n.$
It
is well-known that, for non--dissipative interactions, i.e. when
$\ep=1$, the gain part $\Q^+$ can be written as an integral operator
with explicit kernel \cite{carleman, Mo}
(see also \cite{Grad, cerci} for similar results for the \textit{linearized} Boltzmann equation). We prove that such a representation is still
valid in the dissipative case:
\begin{theo}\label{first} If $f \geq 0$ is such that $\sigma(\v)f(\v) \in L^1(\R,\d\v)$,  then
$$\Q^+(f)(\v)=\int_{\R}k(\v,\v')f(\v')\d\v'$$
where the integral kernel $k(\v,\v')$ can be made explicit (see
\eqref{k}).
\end{theo}

Actually, most important is the fact that the integral
kernel $k(\v,\v')$ turns out to be very similar to that obtained in the classical
case (see for instance \cite{Mo, carleman}), the only changes
standing in some explicit numerical constants. Moreover, as we shall
see, the kernel $k(\v,\v')$ and the Maxwellian distribution
$\mathcal{M}$ satisfy the following {\bf detailed balance law}:
$$k(\v,\v')\mathcal{M}(\v')=k(\v',\v)\mathcal{M}(\v), \qquad \v, \v' \in \R,$$
that allows us to recover Theorem \ref{existence} in a direct way.
Recall  that, in \cite{loto}, the Gaussian nature of the steady
state of $\Q$ was obtained by replacing $\Q$ by its grazing
collision limit.

 We derive from
these two results some important consequences on the linear
Boltzmann equation \eqref{BoltzNL} with $\mathcal{C}=\Q$. The
applications are dealing with the space dependent version of
\eqref{BoltzNL} as well as with the space homogeneous version of it.
The first one concerns the spectral properties of the Boltzmann
collision operator in its natural Hilbert space setting.

\subsection{Spectral properties of the Boltzmann operator in $L^2(\mathcal{M}^{-1})$}
Applying the above $H$--Theo-\\ rem \ref{theoH} with the quadratic
convex function $\Phi(x)=(x-1)^2$, one sees that a natural function
space for the study of the {\it homogeneous linear Boltzmann
equation} is the weighted space $L^2(\mathcal{M}^{-1})$. Now, from
Theorem \ref{first}, it is possible to prove that the gain collision
operator $Q^+$ is {\it compact} in $L^2(\mathcal{M}^{-1})$. This
compactness result has important consequences on the structure of
the spectrum of $\Q$ as an operator in $L^2(\mathcal{M}^{-1})$.
Precisely, from Weyl's Theorem,  the spectrum of $\Q$ in this space
is given by the (essential) range of the collision frequency
$\sigma(\cdot)$  and of isolated eigenvalues with finite algebraic
multiplicities. Since $\lambda=0$ is a simple eigenvalue of $\Q$
(its associated null space is spanned by $\mathcal{M}$), this leads
to the existence of a positive spectral gap. In turns, one proves
that any solution to the space-homogoneous linear Boltzmann equation
\eqref{BoltzNL} converges at an exponential rate towards the
equilibrium. \textit{These spectral results are  technical
generalizations of some of the fundamental results of T. Carleman
\cite{carleman}, but are new in the context of granular gas
dynamics}.

\subsection{Honest solutions for hard--spheres model}
It is easily seen that, for any nonnegative $f$,
\begin{equation}\label{conserva}
\int_{\R}\Q^+(f)(\v)\d \v=\int_{\R}\sigma(\v)f(\v)\d \v,\end{equation}
i.e. the collision operator $\Q$ is \textit{conservative.} Then,
formally, any  nonnegative solution $f(\x,\v,t)$ to \eqref{BoltzNL}
(with $\mathcal{C}=\Q$) should satisfy the following \textit{mass
conservation equation}:
\begin{equation}\label{massconser}
\int_{\R \times \R} f(\x,\v,t)\d \x\d\v=\int_{\R \times \R}
f(\x,\v,0)\d \x \d \v, \qquad \forall t > 0.\end{equation} It is the
main concern of Section \ref{sec:honest} to prove that such a formal
mass conservation property holds true for any nonnegative initial  datum
$f(\x,\v,0) \in  L^1(\R \times \R).$ As well documented in the
monograph \cite{arloban}, this is strongly related to the honesty of
the $C_0$-semigroup governing Eq. \eqref{BoltzNL}. More precisely,
if we denote by $T_0$ the streaming operator:
$$\D(T_0)=\{f \in X\,,\v \cdot \nabla_\x f \in X\}, \qquad T_0f=-\v
\cdot \nabla_\x f,$$ it is not difficult to see that there exists
some extension $G$ of $T_0+\Q$ that generates a $C_0$-semigroup of
contractions $\zt$ in $X=L^1(\R \times \R).$ According to the
so--called "sub-stochastic perturbation" theory, developed in
\cite{Ar, arloban,voigt}, it can be proved that
$$\int_{\R \times \R}Z(t)f(\x,\v)\d \x\d\v=\int_{\R \times \R}f(\x,\v)\d\x\d\v,
\qquad \forall f \in X,\:f \geq 0$$ if and only if $G$ is the
closure of the full transport operator: $G=\overline{T_0+\Q}.$ We
show in Section 4 that the latter holds. To do so, we shall use the
integral representation (Theorem \ref{first}) in order to apply some
of the results of \cite{Ar} (see also \cite[Chapter 10]{arloban})
dealing with the
classical linear Boltzmann equation.  

\subsection{Organization of the paper} We derive in Section 2 the integral representation of $\Q^+$ (Theorem \ref{integraHS})
 as well as
some of its immediate consequences concerning the explicit
expression of the collision frequency. We also recover Theorem
\ref{existence} through a detailed balance law. Section 3 is devoted
to the study of the collision operator $\Q$ in the narrow space
$L^2(\mathcal{M}^{-1}(\v)\d\v)$ and its spectral consequences. In
Section 4 we apply the results of Section 2 as well as some known
facts about the classical linear Boltzmann equation
\cite{Ar,arloban} to  the honesty of the solutions to the Boltzmann
equation for dissipative hard-spheres.

\section{Integral representation of the gain operator}

Let us consider the gain
operator for dissipative  hard--spheres:
$$\Q^+(f)(\v)=\epsilon^{-2}\It |\q \cdot \n|f(\v_{\star})\mathcal{M}_1(\w_{\star})\d \w \d
\n$$ and let $\sigma(\v)$ be the corresponding collision frequency:
$$\sigma(\v)=\It |\q \cdot \n|\mathcal{M}_1(\w)\d \w\d \n,\qquad \v \in \R.$$
Recall  that $\mathcal{M}_1$ is a Maxwellian distribution
function with bulk velocity $\u$ and effective temperature
$\vartheta_1$.
 We recall here the general microscopic description of the pre-collisional
 velocities $(\vb,\wb)$ which result in $(\v,\w)$ after collision. For a
 constant restitution coefficient $0<\epsilon <
1$, one has \cite{poschel,villani}
 \begin{equation*}
\begin{cases}
\v_{\star}=\v-2\alpha\dfrac{1-\beta}{1-2\beta} [\q \cdot \n] \n,\\
\w_{\star}=\w+2(1-\alpha)\dfrac{1-\beta}{1-2\beta} [\q \cdot \n] \n;
\end{cases}\end{equation*}
where $\q=\v-\w$, $\alpha$ is the mass ratio and $\beta$ denotes the
inelasticity parameter
$$\alpha=\dfrac{m_1}{m+m_1}, \qquad \qquad \qquad
\beta=\dfrac{1-\epsilon}{2}.$$

We show in this section that, as it occurs for the classical
Boltzmann equation, $\Q^+$ turns out to be an integral operator with
explicit kernel. The proof of such a result is based on well-known
tools from the linear elastic scattering theory \cite{carleman,
Grad, Mo} while, in the dissipative case, similar calculations have
been performed to derive a Carleman representation of the nonlinear
Boltzmann operator in \cite{mischler}.

\begin{theo}[\textbf{Integral representation of}
$\Q^+$]\label{integraHS} For any $f \in L^1(\R \times \R,\d x
\otimes \sigma(\v)\d \v)$, \begin{equation} \label{Q+}
\Q^+f(\x,\v)=\dfrac{1}{2\epsilon^2\gamma^2}\int_{\R}f(\x,\v')k(\v,\v')\d
\v',\end{equation} where \begin{equation}\label{k}
k(\v,\v')=\left(\dfrac{m_1}{2\pi
\vartheta_1}\right)^{1/2}|\v-\v'|^{-1}\exp
\left\{-\dfrac{m_1}{8\vartheta_1} \left(
(1+\mu)|\v-\v'|+\dfrac{|\v-\u|^2-{|\v'-\u|}^2}{|\v-\v'|}\right)^2\right\}\end{equation}
with $\mu=-\dfrac{2\alpha(1-\b)-1}{\alpha(1-\b)}>0$ and
$\gamma=\alpha \dfrac{1-\b}{1-2\b}.$
\end{theo}
\begin{proof} The local (in $\x$) nature of $\Q^+$ is obvious and we can restrict ourselves to prove the result
for a function  $f \in L^1(\R,\sigma(\v)\d\v)$ that does not depend
on $\x$. Set $\gamma=\alpha \frac{1-\b}{1-2\b}$ and
$\g=(1-\alpha)\frac{1-\b}{1-2\b}$ so that
$$\vb=\v-2\gamma [\q \cdot \n]\n \quad \text{ and } \quad \wb=\w+2\g\, [\q \cdot
\n]\n.$$
The following formula, for smooth $\varphi$:
\begin{equation*}\label{formule1}
\int_{\mathbb{S}_+}(\q \cdot \n)\varphi\left((\q \cdot
\n)\n\right)\d
\n=\dfrac{|\q|}{4}\int_{\S}\varphi\left(\dfrac{\q-|\q|\sigma}{2}\right)\d
\sigma=\dfrac{1}{2}\int_{\R}\delta(2 \x \cdot \q +
\x^2)\varphi(-\x/2)\d \x.\end{equation*} applied to
$$\varphi(x)=f\left(v-2\gamma x\right)\mathcal{M}_1
\left(\w +2\overline{\gamma} x\right)$$ yields
%
\begin{equation*}
\Q^+f(\v)=\epsilon^{-2}\IRR \delta(2 \x \cdot \q + \x^2)f(\v+\gamma
\x)\mathcal{M}_1(\w-\g \x)\d \w \d \x.\end{equation*}
The change of variables $\x \mapsto \v'=\v + \gamma \x$ leads to
$$\Q^+f(\v)=\epsilon^{-2}\gamma^{-3}\IRR \delta(2
\gamma^{-1}(\v'-\v) \cdot \q +
\gamma^{-2}|\v'-\v|^2)f(\v')\mathcal{M}_1(\w-\dfrac{\g}{\gamma}(\v'-\v))\d
\w \d \v'.$$ Now, keeping $\v$ and $\v'$ fixed, we perform the
change of variables $\w \mapsto \w'=\w-\dfrac{\g}{\gamma}(\v'-\v)$,
which leads to
\begin{multline*}
\Q^+f(\v)=\epsilon^{-2}\gamma^{-3}\IRR \delta\left(2
\gamma^{-1}(\v'-\v) \cdot [\v-\w'-\dfrac{\g}{\gamma}(\v'-\v)] +
\gamma^{-2}|\v'-\v|^2\right)\times \\
\times f(\v')\mathcal{M}_1(\w')\d \w' \d \v'.\end{multline*} Writing
$\w'=\v+\l_1\n+V_2$ with $\lambda_1=(\w'-\v) \cdot \n \in
\mathbb{R},$ $\n=(\v'-\v)/|\v'-\v|$ and $V_2 \cdot \n=0$,
we get, noting that $\d \w'=\d V_2 \d \l_1$,
\begin{multline*}
\Q^+f(\v)
=\epsilon^{-2}\gamma^{-3}\int_{\R}f(\v')\d \v' \int_{\mathbb{R}}
\d\l_1 \int_{V_2 \cdot
n=0}\mathcal{M}_1(\v+V_2+\l_1\n)\d
V_2\times\\
\times \delta\left(\gamma^{-2}|\v'-\v|^2-2 \g \gamma^{-2}|\v'-\v|^2
-2\gamma^{-1}\l_1|\v'-\v|\right).
\end{multline*}
Thanks to the change of variables $\l_1 \mapsto
2\gamma^{-1}|\v'-\v|\l_1$, one can evaluate the Dirac mass as
\begin{multline*}
\int_{\mathbb{R}}\delta\left(\gamma^{-2}|\v'-\v|^2-2 \g
\gamma^{-2}|\v'-\v|^2
-2\gamma^{-1}\l_1|\v'-\v|\right)\mathcal{M}_1(\v+V_2+\l_1\n)\d
\l_1\\
=\dfrac{\gamma}{2|\v'-\v|}\mathcal{M}_1\left(\v+V_2+\dfrac{1-2\g}{2\gamma}(\v'-\v)\right)\end{multline*}
where we used that $\n=(\v'-\v)/|\v'-\v|$. Consequently,
\begin{equation*}
\Q^+f(\v)=\dfrac{1}{2\epsilon^2\gamma^2}\int_{\R}k(\v,\v')f(\v')\d
\v'\end{equation*} where
$$k(\v,\v')=\dfrac{1}{|\v-\v'|}\int_{V_2 \cdot
(\v'-\v)=0}\mathcal{M}_1\left(\v+V_2+\dfrac{1-2\g}{2\gamma}(\v'-\v)\right)\d
V_2.$$ It remains now to explicit $k(\v,\v')$. We will use
the approach of \cite{Mo}. Let us assume $\v,\v'$ to be fixed. Let
$\mathsf{P}$ be the hyperplan orthogonal to $(\v'-\v)$. For any
$V_2 \in \mathsf{P}$, set
$$z=\v+\dfrac{1-2\g}{2\gamma}(\v'-\v)+V_2-\u$$
so that
$$k(\v,\v')=\left(\frac{\varrho_1}{\pi}\right)^{3/2}|\v-\v'|^{-1}\int_{V_2 \in \mathsf{P}}\exp\{-\varrho_1 z^2\}\d V_2.$$
where $\varrho_1=\dfrac{m_1}{2\vartheta_1}.$ Denoting for simplicity
$u=\dfrac{\v+\v'}{2}-\u$ and $\mu=-\dfrac{1-2\g}{\gamma}$, one has
\begin{equation*}\begin{split}
z^2&=\left(u+\dfrac{\v-\v'}{2}+\dfrac{\mu}{2}(\v-\v')+V_2\right)^2\\
&=|u+V_2|^2+\dfrac{(1+\mu)^2}{4} |\v-\v'|^2 +
\dfrac{1+\mu}{2}(|\v-\u|^2-|\v'-\u|^2) \end{split}\end{equation*} where we
used the fact that $V_2$ is orthogonal to $(\v'-\v)$. Splitting $u$ as
$$u=u_0+u_\bot$$
where $u_0$ is parallel to $\v-\v'$ while $u_\bot$
is orthogonal to $\v-\v'$ (i.e. $u_\bot \in \mathsf{P}$),
we see that
$$|u+V_2|^2=|u_0|^2+|u_\bot+V_2|^2
\qquad \text{ and } \qquad |u_
0|^2=\dfrac{\left[|\v-\u|^2-{|\v'-\u|}^2\right]^2}{4|\v-\v'|^2},
$$
so that
\begin{multline*}
k(\v,\v')=|\v-\v'|^{-1}\left(\frac{\varrho_1}{\pi}\right)^{3/2}
\int_{\mathsf{P}}\exp\left(-\varrho_1|u_\bot+V_2|^2\right)\d
V_2\\
\exp \left\{-\dfrac{\varrho_1}{4} \left(
{(1+\mu)^2}|\v-\v'|^2+2(1+\mu)(|\v-\u|^2-|\v'-\u|^2)+\dfrac{\left[|\v-\u|^2-{|\v'-\u|}^2\right]^2}{|\v-\v'|^2}\right)\right\}
.\end{multline*} Finally, since $u_\bot \in
\mathsf{P}$,
$$\int_{\mathsf{P}}\exp\left(-\varrho_1|u_\bot+V_2|^2\right)\d
V_2=\int_{\mathbb{R}^2}\exp(-\varrho_1 \x^2)\d
\x=\dfrac{\pi}{\varrho_1},$$ one obtains the desired expression for
$k(\v,\v')$.
\end{proof}


The very important fact to be noticed out is that the expression of
$k(\v,\v')$ is very similar to that one obtains in the elastic case
\cite{Mo}, the only change being the expression of the constant
$\mu$. In particular, in the elastic case $\epsilon=1$, we recover
the expression of the kernel obtained in \cite{carleman} for
particles of same mass (i.e. $m=m_1$) and in \cite{Mo} for particles
with different masses.

Another fundamental property of the kernel $k(\v,\v')$
is that it allows us to recover the steady state of $\Q$ through some {\it microscopic
  detailed balance law}. Precisely,
\begin{theo}\label{theo:balance} With the notations of the Theorem \ref{integraHS}, the following
  \textbf{\textit{detailed balance law}}:
\begin{equation}\label{balance}
k(\v,v')\exp\left\{-\dfrac{m_1}{2\vartheta_1}(1+\mu)(\v'-\u)^2\right\}=k(\v',\v)\exp\left\{-\dfrac{m_1}{2\vartheta_1}(1+\mu)(\v-\u)^2\right\},
\end{equation}
holds for any $\v,\v' \in \R.$
As  a consequence, the Maxwellian velocity distribution:
$$\mathcal{M}(\v)=\left(\dfrac{m}{2\pi
\vartheta^{\#}}\right)^{3/2}\exp
\left\{-\dfrac{m(\v-\u)^2}{2\vartheta^{\#}}\right\} \qquad \v \in
\R,$$ 
with $ \vartheta^{\#}=\dfrac{(1-\a)(1-\b) }{1-\a(1-\b)}\vartheta_1$
is the unique  equilibrium state of $\Q$ with unit mass.
\end{theo}
\begin{proof} According to Eq. \eqref{k}, it is easily seen that
$$k(\v',\v)=k(\v,v')\exp\left\{
 \frac{m_1}{2\vartheta_1}(1+\mu)\left(|\v-\u|^2-|\v'-\u|^2\right)\right\},\qquad \v,\v' \in \R$$
which is nothing but  \eqref{balance}. Now, writing
$\frac{m_1}{2\vartheta_1}(1+\mu)=\frac{m}{2\vartheta^\sharp}$,
straightforward calculations lead to the desired expression for the
equilibrium temperature $\vartheta^\sharp.$ The  fact that
$\mathcal{M}$ is an equilibrium solution with unit mass follows then
from the fact that
$$\Q(\mathcal{M})(\v)=\int_{\R}k(\v,\v')\mathcal{M}(\v')\d\v'-\sigma(\v)\mathcal{M}(\v)=\int_{\R}\left[k(\v,\v')\mathcal{M}(\v')-k(\v',\v)\mathcal{M}(\v)\right]\d\v'$$
and from the detailed balance law \eqref{balance}. To prove that the steady state is unique, we adopt the stategy of \cite[Theorem 1]{poupaud}. Precisely, consider the equation
\begin{equation}\label{equa:poupaud}
\sigma(\v)f(\v)=\Q^+f(\v),\qquad \forall \v \in \R
\end{equation}
which admits at least the solution $f=\mathcal{M}.$ Since
$\sigma(\v)$ does not vanish, any solution $f$ to
\eqref{equa:poupaud} is such that
$$f(\v)=\frac{1}{\sigma(\v)}\,\Q^+(f)(\v),\qquad \qquad \forall \v \in \R.$$
Since $Q^+$ is an integral operator with \textit{nonnegative}
kernel, it is clear that $\sigma(\v)|f(\v)| \leq \Q^+(|f|)(\v)$ for
any $\v \in \R$. Now, from the positivity of both $\sigma$ and
$\Q^+$,  one sees that the conservation of mass \eqref{conserva}
reads:
\begin{equation*}
\|\sigma f\|_{X} =\|\sigma\,|f|\,\|_{X}=\|\Q^+(|f|)\|_{X}.
\end{equation*}
This shows that, actually,
$|\Q^+(f)(\v)|=\sigma(\v)|f(\v)|=\Q^+(|f|)(\v)$ for any $\v \in \R.$
Again, since $\Q^+$ is a positive operator, one obtains that
$$f=\pm|f|.$$
Now, assume that \eqref{equa:poupaud} admits two solutions $f_1$,
$f_2$ with $\int_{\R} f_1(\v)\d\v=\int_{\R}f_2(\v)\d\v=1.$ Then,
$f_1-f_2$ is again a solution to \eqref{equa:poupaud} so that,
$f_1-f_2=\pm |f_1-f_2|.$ Thus,
$$\pm \int_{\R}|f_1(\v)-f_2(\v)|\d\v=\int_{\R} f_1(\v)\d\v-\int_{\R}f_2(\v)\d\v=0$$
and the uniqueness follows.
\end{proof}

The above result allows to derive the explicit expression of the collision frequency $\sigma(\v)$:
\begin{cor} The collision frequency $\sigma(\v)$ for dissipative
  hard--spheres interactions is given by
\begin{multline}
\label{sigmaHS} \sigma(\v)=\dfrac{2\pi}{(2+\mu)^2}
\sqrt{\dfrac{m_1}{2\pi\vartheta_1}}\left\{\frac{4\vartheta_1}{m_1}
\exp\left(-\frac{m_1}{2\vartheta_1}|\v-\u|^2\right) \right. \\
\left. +\left(2|\v-\u|+\frac{2\vartheta_1}{m_1|\v-\u|}\right)
\int_0^{2|\v-\u|}\exp\left(-\frac{m_1}{8\vartheta_1}t^2\right)\d t
\right\}.
\end{multline}
Consequently, there exist positive constants $\nu_0,\nu_1$ such that
$$\nu_0(1+|\v-\u|) \leq \sigma(\v) \leq \nu_1 (1+|\v-\u|), \qquad \forall
\v \in \R.$$\end{cor}
\begin{proof} Set $C=\sqrt{\dfrac{m_1}{2\pi\vartheta_1}}$.  Noting
that $\sigma(\v)=\int_{\R}k(\v',\v)\d\v'$ for any $\v \in \R,$ one
has, with the change of variable $z=\v'-\v$, in a polar coordinate
system in which $\v$ lies on the third axis
\begin{equation*}\begin{split}
\sigma(\v)&=C\int_{\R}\exp\left\{-\dfrac{m_1}{8\vartheta_1}\left((1+\mu)|z|
- \dfrac{|\v-\u|^2-|z+\v-\u|^2}{|z|}\right)^2\right\}|z|^{-1}\d z\\
&=2\pi C\int_0^\infty \d \varrho \int_0^{\pi} \varrho
\exp\left\{-\dfrac{m_1}{8\vartheta_1}\left((2+\mu)\varrho+2|\v-\u|\cos
\varphi\right)^2\right\}\sin \varphi \,\d\varphi.
\end{split}\end{equation*}
The computation of this last integral leads to the desired
expression for $\sigma(\v).$ The estimates are then straightforward \cite{Mo}.
\end{proof}

\section{Application to the Boltzmann operator in
$L^2(\mathcal{M}^{-1}).$}

We investigate in this section the properties of the Boltzmann
operator $\Q$ in the weighted space
$$\mathcal{H}=L^2(\R;\mathcal{M}^{-1}(\v)\d\v).$$
We shall denote by $\langle \cdot,\cdot \rangle_{\mathcal{H}}$ the
inner product in $\mathcal{H}$. The introduction of such an Hilbert
space setting is motivated by the application of the $H$-Theorem
\ref{theoH} with the convex function $\Phi(x)=(x-1)^2.$ In this
case, one sees that, if $f_0 \geq 0$ is a space homogeneous initial
distribution  such that
$$\int_{\R}f_0(\v)\d\v=1,\qquad
\int_{\R}|f_0(\v)|^2\mathcal{M}^{-1}(\v)\d\v < \infty,$$ then any solution $f(t,\v)$ to the space homogeneous equation
\begin{equation}\label{spacehomo}
\partial_t f(t,\v)=\Q(f)(t,\v),\qquad f(0,\v)=f_0(\v) \in
\mathcal{H},\end{equation} satisfies the following estimate:
$$\dfrac{\d}{\d t}\int_{\R}\left|f(t,\v)-\mathcal{M}(\v)\right|^2\mathcal{M}(\v)^{-1}\d\v
\leq 0,\qquad \qquad    t \geq 0.$$ In other words, the mapping
$t \longmapsto \left\|f(t,\cdot)-\mathcal{M}\right\|_{\mathcal{H}}$
is {\it nonincreasing}. For these reasons, the study of the properties of the collision operator $\Q$ in
$\mathcal{H}$ is of particular relevance for the asymptotic behavior of the solution
\begin{equation}\label{homogeneous}
\partial_t f(t,\v)=\Q(f)(t,\v),\qquad f(0,\v) \in
\mathcal{H}.\end{equation}
The material of this section borrows some techniques already employed by T. Carleman  \cite{carleman} in the study of
non--dissipative gas dynamics (see also, e.g. \cite{cerci} or \cite{Grad} for similar results
in the context of the \textit{linearized} Boltzmann equation).
Let  $\mathcal{L}$ be the realization of the operator $\Q$ in $\mathcal{H}$, i.e.
\begin{equation*} \D(\mathcal{L})=\left\{f \in \mathcal{H}\,;\int_{R}|f(\v)|^2\sigma(\v)\mathcal{M}^{-1}(\v)\d\v
<\infty\right\}.\end{equation*}
and, for any $f \in \D(\mathcal{L})$, $\mathcal{L}f(\v)=\Q(f)(\v)$ is given by \eqref{linearoperator}. As previously, one can use
the following splitting of $\mathcal{L}$ as a gain operator and a
loss (multiplication) operator,
$\mathcal{L}=\mathcal{L}^+-\mathcal{L}^-$ with
$$\mathcal{L}^+(f)(\v)=\int_{\R}k(\v,\v')f(\v')\d\v'\qquad \text{ and } \qquad \mathcal{L}^-(f)=\sigma(\v)f(\v), \qquad \qquad f \in \D(\mathcal{L}).$$
We shall show, as in the classical case, that $\mathcal{L}^+$ is
actually a bounded operator in $\mathcal{H}.$ Precisely, let
$\mathscr{J}$ define the natural bijection operator from
$L^2(\R,\d\v)$ to $\mathcal{H}$:
\begin{equation*}\begin{cases}
\mathscr{J}\::\:&L^2(\R,\d\v)                  \longrightarrow \mathcal{H}\\
&f \longmapsto \mathscr{J}f(\v)=\mathcal{M}^{1/2}(\v)f(\v)
                 \end{cases}
\end{equation*}
It is clear that $\mathscr{J}$ is a bounded bijective operator whose
inverse is given by
$$\mathscr{J}^{-1}g(\v)=\mathcal{M}^{-1/2}(\v)g(\v) \in L^2(\R,\d\v), \qquad \forall g \in \mathcal{H}.$$
Now, let us define
$$G(\v,\v')=\mathcal{M}^{-1/2}(\v)k(\v,\v')\mathcal{M}^{1/2}(\v'),\qquad \qquad \v,\v' \in \R,$$
i.e.
\begin{equation}\label{gvv}
G(\v,\v')=\left(\dfrac{m_1}{2\pi
\vartheta_1}\right)^{1/2}|\v-\v'|^{-1}\exp
\left\{-\dfrac{m_1}{8\vartheta_1} \left(
(1+\mu)^2|\v-\v'|^2+\dfrac{(|\v-\u|^2-{|\v'-\u|}^2)^2}{|\v-\v'|^2}\right)\right\}.
\end{equation}
From the detailed balance law \eqref{balance}, one easily checks
that $G(\v,\v')=G(\v',\v)$ for any $\v,\v' \in \R$. Therefore,
defining $\mathcal{G}$ as the  integral operator
in $L^2(\R,\d\v)$ with kernel $G(\v,\v')$, i.e.
$$\mathcal{G}f(v)=\int_{\R}G(\v,\v')f(\v')\d\v',$$
one can prove the following:
\begin{propo}\label{prop:l+} $\mathcal{G}$ is a bounded symmetric operator in $L^2(\R,\d\v)$
and $\mathcal{L}^+=\mathscr{J}\mathcal{G}\mathscr{J}^{-1}.$
 Consequently, $\mathcal{L}^+$ is a bounded symmetric operator in $\mathcal{H}$.
\end{propo}
\begin{proof} It is clear that $\mathcal{G}$ is symmetric since $G(\v,\v')=G(\v',\v)$. Now, to prove
the boundedness of $\mathcal{G}$, one adopts a strategy already used
in the non--dissipative case by T. Carleman \cite[p. 75]{carleman}
and shows easily that
$$C:=\sup_{\v \in \R}\int_{\R}
G(\v,v')\d\v' < \infty.$$ Since $G(\cdot,\cdot)$ is symmetric, one
also has $\sup_{\v'\in \R}\int_{\R}G(\v,\v')\d\v=C <\infty$.
Denoting by $\langle \cdot,\cdot \rangle$ the usual inner product of
$L^2(\R,\d\v)$, one deduces from Cauchy-Schwarz identity,
$$\langle \mathcal{G}f,g\rangle \leq \dfrac{C}{2}\left(\int_{\R}|f(\v)|^2\d\v + \int_{\R}|g(\v')|^2\d\v'\right),\qquad \qquad \forall f,g \in L^2(\R,\d\v),$$
which leads to the boundedness of $\mathcal{G}.$ Now, since
$G(\v,\v')=\mathcal{M}^{-1/2}(\v)k(\v,\v')\mathcal{M}^{1/2}(\v')$ for any
$\v,\v' \in \R,$ one gets easily that
$\mathcal{L}^+=\mathscr{J}\mathcal{G}\mathscr{J}^{-1}$ and the
conclusion follows.\end{proof}

In Proposition \ref{prop:l+}, we proved that the gain operator
$\mathcal{L}^+$ is bounded in $\mathcal{H}$, i.e. $\mathcal{L}^+ \in
\mathscr{B}(\mathcal{H})$. Actually, we have much better and it is
possible, as in the non--dissipative case, to prove that
$\mathcal{L}^+$ is a compact operator in $\mathcal{H}$. Precisely,
the following lemma is a direct consequence of Theorem
\ref{integraHS} and similar calculations valid for the
non-dissipative case \cite[p. 70--75]{carleman}. However, we give a
detailed proof of it since the known similar results  by  T.
Carleman  are all dealing with the case $m=m_1$ and $\epsilon=1$. It
has to be checked that taking account the parameters $m \neq m_1$
and $\epsilon < 1$ does not lead to supplementary difficulty (see
Remark \ref{eigen} where the role of $\epsilon \neq 1$ does not
allow to adapt {\it mutatis mutandis} a result valid in the elastic
case).
\begin{lemme}
For any $0 < p < 3$ and any $q \geq 0$, there exists $C(p,q) >0$
such that
$$\int_{\R}|G(\v,\v')|^p \dfrac{\d \v'}{(1+|\v'-\u|)^q} \leq \dfrac{C(p,q)}{(1+|\v-\u|)^{q+1}}, \qquad \forall \v \in \R.$$
\end{lemme}
\begin{proof}
The proof is a technical generalization of a similar result due to
T. Carleman \cite{carleman} in the classical case (i.e. when $m=m_1$
and $\epsilon=1$).  Let us fix $0 < p < 3$ and $q \geq 0$ and set
$$I(\v)=\int_{\R}|G(\v,\v')|^p \dfrac{\d \v'}{(1+|\v'-\u|)^q}.$$
Then, one sees easily that
$$I(\v)=
2\pi\left(\dfrac{m_1}{2\pi \vartheta_1}\right)^{p/2} \int_0^\pi \sin
\varphi \d\varphi \int_0^\infty
\varrho^{2-p}\dfrac{\exp\left\{-\dfrac{m_1p}{8\vartheta_1}
\left((1+\mu)^2\varrho^2+(\varrho+2|\v-\u|\cos
\varphi)^2\right)\right\}}{\left(1+\sqrt{\varrho^2+|\v-\u|^2+2\varrho|\v-\u|\cos
\varphi}\right)^q}\d\varrho.$$
Note that, since $0< p < 3$,
\begin{equation}\label{estim1}
 \sup_{\v \in \R} I(\v) \leq 4\pi\left(\dfrac{m_1}{2\pi
\vartheta_1}\right)^{p/2} \int_0^\infty
\varrho^{2-p}\exp\left\{-\dfrac{m_1p}{8\vartheta_1}
(1+\mu)^2\varrho^2\right\}\d\varrho<\infty.
\end{equation}
Performing the change of variable $x=\varrho/|\v-\u| + 2\cos
\varphi,$ $y=\varrho/|\v-\u|$, one has $(x,y) \in \Omega$ where
$$\Omega=\{(x,y) \in \mathbb{R}^2\,;\,y >0,\,|x-y|\leq 2\}$$
and
$$I(\v)=\left(\dfrac{m_1}{2\pi
\vartheta_1}\right)^{p/2}\pi|\v-\u|^{3-p}\int_{\Omega}
\exp\left\{-\dfrac{m_1p|\v-\u|^2}{8\vartheta_1} \left((1+\mu)^2
y^2+x^2\right)\right\}\dfrac{\d x \d
y}{y^{p-2}\left(1+|\v-\u|\sqrt{1+xy}\right)^q}.$$ We split $\Omega$
into $\Omega=\Omega_1 \cup \Omega_2$ where $\Omega_1$ is the
half-ellipse $$\Omega_1=\{(\x,y) \in \mathbb{R}^2\,:\,y >
0,\,(1+\mu)^2y^2+x^2 < 1/4\}\qquad \text{ while }\qquad
\Omega_2=\Omega \setminus \Omega_1.$$ Note that, since $1+\mu \geq
1,$ one has $\Omega_1 \subset \Omega.$ One defines correspondingly
$I_1(\v)$ and $I_2(\v)$ as the above integral over $\Omega_1$ and
$\Omega_2$ respectively. One notes first that, if $(x,y) \in
\Omega_1$ then $xy > -\frac{1}{8(1+\mu)}$ so that
$$I_1(\v) \leq  \left(\dfrac{m_1}{2\pi
\vartheta_1}\right)^{p/2}\dfrac{\pi|\v-\u|^{3-p}}{(1+a|\v-\u|)^q}\int_{\Omega_1}
\exp\left\{-\dfrac{m_1p|\v-\u|^2}{8\vartheta_1} \left((1+\mu)^2
y^2+x^2\right)\right\}\dfrac{\d x \d y}{y^{p-2}}$$ where
$a=\sqrt{1-\frac{1}{8(1+\mu)}}$, $0<a<1$. Letting $R=\left(\frac{m_1
p}{8 \vartheta_1}\right)^{1/2}$ and setting $t=R|\v-\u|x,$
$u=R(1+\mu)|\v-\u|y$, it is easy to check that $-R|\v-\u|/2 \leq t
\leq R|\v-\u|/2,$ while $0 \leq u \leq R|\v-\u|/2,$ so that
$$I_1(\v) \leq \left(\dfrac{m_1}{2\pi
\vartheta_1}\right)^{p/2} \dfrac{\pi
R^{p-4}(1+\mu)^{p-3}}{|\v-\u|(1+a|\v-\u|)^q}\int_{\mathbb{R}}\d t
\int_0^{\infty} \dfrac{\exp\{-(t^2+u^2)\}}{u^{p-2}}\d u.$$ Thus,
there exists a constant $C_1(p,q) >0$ such that
\begin{equation}\label{estim11}
I_1(\v) \leq \dfrac{C_1(p,q)}{|\v-\u|(1+a|\v-\u|)^q}, \qquad \forall
\v \in \R.
\end{equation}
Let us now deal  with $I_2(\v)$. Arguing as above,
\begin{equation*}\begin{split}
I_2(\v)
&=\left(\frac{m_1}{2\pi
\vartheta_1}\right)^{p/2}\pi|\v-\u|^{3-p}\int_{\Omega_2}
\exp\left\{-\frac{R^2|\v-\u|^2}{2} \left((1+\mu)^2
y^2+x^2\right)\right\}\times \\
&\phantom{++++ ++++++}\times
\frac{\exp\left\{-\frac{R^2|\v-\u|^2}{2} \left((1+\mu)^2
y^2+x^2\right)\right\}}{y^{p-2}\left(1+|\v-\u|\sqrt{1+xy}\right)^q}\d
x \d y.
\end{split}\end{equation*}
Clearly, since $(1+\mu)^2y^2+x^2 >1/4$ for any $(x,y) \in \Omega_2$,
then
\begin{equation*}\begin{split}
I_2(\v)&\leq \left(\dfrac{m_1}{2\pi
\vartheta_1}\right)^{p/2}\pi|\v-\u|^{3-p}\int_{\Omega_2}
\exp\left(-\frac{R^2|\v-\u|^2}{8}\right)\exp\left\{-\frac{R^2|\v-\u|^2}{2}
\left((1+\mu)^2 y^2+x^2\right)\right\}\dfrac{\d x \d
y}{y^{p-2} }\\
&\leq \left(\dfrac{m_1}{2\pi \vartheta_1}\right)^{p/2} \pi
|\v-\u|^{3-p}\exp\left(-\frac{R^2|\v-\u|^2}{8}\right)\int_0^\infty
\exp\left(-\frac{R^2|\v-\u|^2}{2}(1+\mu)^2 y^2  \right) \dfrac{\d
y}{y^{p-2}}\int_{y-2}^{y+2}\d x.
\end{split}\end{equation*}
Hence, there is some constant $C_2(p,q)$ such that
\begin{equation}\label{i2}
I_2(\v) \leq C_2(p,q) \exp\left(-\frac{R^2|\v-\u|^2}{8}\right),
\qquad \v \in \R.\end{equation} Combining \eqref{estim11} and
\eqref{i2}, one sees that
$$I(\v) \leq \dfrac{C_1(p,q)}{|\v-\u|(1+a|\v-\u|)^q} +C_2(p,q)\exp\left(-\frac{R^2|\v-\u|^2}{8}\right), \qquad \v \in
\R.$$ According to \eqref{estim1},   $\limsup_{|\v-\u| \to 0}I(\v) <
\infty,$ from which we get the conclusion.
\end{proof}

\begin{nb}
Note that the above Lemma can be extended to more general collision
kernels (including long-range interactions) following the lines of
the recent results \cite{strain} dealing with the elastic case.
\end{nb}

From the above Lemma, one has the following compactness result:
\begin{propo}\label{propo:compact} $\mathcal{G}$ is compact in $L^2(\R,\d\v)$. Consequently, $\mathcal{L}^+$ is a compact operator in $\mathcal{H}.$
\end{propo}
\begin{proof}
Applying arguments already used in \cite{carleman}, the above Lemma
implies that the third iterate of $\mathcal{G}$ is an
Hilbert--Schmidt operator in $L^2(\R,\d\v)$, i.e. the kernel of
$\mathcal{G}^3$ is square summable over $\R \times \R.$  The
compactness of $\mathcal{G}$ follows then from standard arguments
and that of $\mathcal{L}^+$ is deduced from the identity
$\mathcal{L}^+=\mathscr{J}\mathcal{G}\mathscr{J}^{-1}$ (see
Proposition \ref{prop:l+}).
\end{proof}
The following, which generalizes a known result from classical kinetic theory,
proves that $\mathcal{L}$ is a negative symmetric operator in $\mathcal{H}$:
\begin{propo}\label{negat}
The operator $(\mathcal{L}, \D(\mathcal{L}))$ is a negative self--adjoint operator of $\mathcal{H}$. Precisely,
$$\langle \mathcal{L}f,f\rangle_{\mathcal{H}}=-\dfrac{1}{2}\int_{\R \times \R}
k(\v,\v')\mathcal{M}(\v')\bigg[\mathcal{M}^{-1}(\v)f(\v)-f(\v')\mathcal{M}^{-1}(\v')\bigg]^2\d\v\d\v'\leq 0$$
 for any $f \in \D(\mathcal{L})$. \end{propo}
\begin{proof} The fact that $(\mathcal{L},\D(\mathcal{L}))$ is self-adjoint is a direct consequence of Proposition
 \ref{prop:l+} since $\mathcal{L}^-$ is clearly symmetric.
Now, it is a classical feature, from the detailed balance law
\eqref{balance}, that
$$\langle \mathcal{L}f,f\rangle_{\mathcal{H}}=\int_{\R \times
\R}k(\v,\v')\mathcal{M}(\v')\bigg[\mathcal{M}^{-1}(\v')f(\v')-f(v)
\mathcal{M}^{-1}(\v)\bigg]f(\v)\mathcal{M}^{-1}(\v)\d\v\d\v'.$$
Exchanging $\v$ and $\v'$ and using again the detailed balance law
\eqref{balance}, one sees that
$$ \langle \mathcal{L}f,f\rangle_{\mathcal{H}}
=\int_{\R \times
\R}k(\v,\v')\mathcal{M}(\v')\bigg[\mathcal{M}^{-1}(\v)f(\v)-f(\v')
\mathcal{M}^{-1}(\v')\bigg]f(\v')\mathcal{M}^{-1}(\v')\d\v\d\v'$$
so that, taking the mean of the two quantities,
$$\langle \mathcal{L}f,f\rangle_{\mathcal{H}}=-\dfrac{1}{2}\int_{\R \times \R}k(\v,\v')
\mathcal{M}(\v')\bigg[\mathcal{M}^{-1}(\v)f(\v)-f(\v')\mathcal{M}^{-1}(\v')\bigg]^2\d\v\d\v'\leq 0$$
which ends the proof.\end{proof}

\begin{nb} From the above result, the spectrum $\mathfrak{S}(\mathcal{L})$ of $\mathcal{L}$ lies in $\mathbb{R}_-$, i.e. $\mathfrak{S}(\mathcal{L}) \subset (-\infty,0]$. It is clear that $\lambda=0$ lies in $\mathfrak{S}(\mathcal{L})$. Precisely $0$ is a
simple eigenvalue of $\mathcal{L}$  since $\mathcal{M}$ is the unique (up to a multiplication factor) steady state of $\mathcal{L}$.\end{nb}
Combining the above results with Proposition \ref{prop:l+} leads to a precise description of the spectrum of $\mathcal{L}$:
\begin{theo}\label{spectral} The spectrum of $\mathcal{L}$ (as an
operator in $\mathcal{H}$) consists of the spectrum of
$-\mathcal{L}^-$ and of, at most, eigenvalues of finite
multiplicities. Precisely, setting $\nu_0=\inf_{\v \in \R}\sigma(\v)
>0$,
$$\mathfrak{S}(\mathcal{L})=\{\lambda \in \mathbb{R}\,;\, \lambda \leq -\nu_0\}
\cup \{\lambda_n\,;\,n \in I\}$$ where $I \subset \mathds{N}$ and
$(\lambda_n)_n$ is a decreasing sequence of real eigenvalues of
$\mathcal{L}$ with finite algebraic multiplicities: $\lambda_0=0 >
\lambda_1 > \lambda_2 \ldots > \lambda_n > \ldots,$
 which unique possible cluster point is  $-\nu_0$.
\end{theo}

\begin{proof}
From Proposition \ref{propo:compact}, $\mathcal{L}$ is nothing but a
compact perturbation of the loss operator $-\mathcal{L}^-$. Hence,
Weyl's Theorem asserts that $\mathfrak{S}(\mathcal{L}) \setminus
\mathfrak{S}(-\mathcal{L}^-)$ consists of, at most, eigenvalues of
finite algebraic multiplicities which unique possible cluster point
is $\sup \{\lambda,\, \lambda \in \mathfrak{S}(-\mathcal{L}^-)\}.$
In particular, up to a rearrangement, one can write
$\mathfrak{S}(\mathcal{L}) \setminus
\mathfrak{S}(-\mathcal{L}^-)=\{\lambda_n,\,n \in I\}$ with
$\lambda_0  > \lambda_1 > \lambda_2 \ldots > \lambda_n \geq \ldots.$
We already saw that $\lambda_0=0$ since $\mathcal{M}$ is a steady
state of $\Q$ and $\mathcal{M} \in \mathcal{H}$. Now, since
$-\mathcal{L}^-$ is a multiplication operator by the collision
frequency $-\sigma(\cdot)$, its spectrum
$\mathfrak{S}(-\mathcal{L}^-)$ is given by the essential range
$R_{\mathrm{ess}}(-\sigma(\cdot))$ of the collision frequency. From
Corollary \ref{sigmaHS}, one sees without difficulty that
$$R_{\mathrm{ess}}(-\sigma(\cdot))=(-\infty,-\nu_0]$$
where $\nu_0=\inf_{\v \in \R}\sigma(\v)=\lim_{|\v-\u| \to 0
}\sigma(\v)=\frac{8}{(2+\mu)^2}\sqrt{\frac{2\pi\vartheta_1}{m_1}}$
is positive.\end{proof}

\begin{nb}\label{eigen} We conjecture that, as it is the case for elastic
interactions \cite{kuscer}, the set of eigenvalues lying in
$(-\nu_0,0)$ is infinite. However, the technical generalization of
the proof of \cite{kuscer} appears to be non trivial because of the
non zero parameter $\mu$. We thank anyway an anonymous referee for
having pointed to us the reference \cite{kuscer}.\end{nb}

The above result provides a complete picture of the spectrum of $\Q$
as an operator in $\mathcal{H}$ (see Fig. 1) and shows, in
particular, the existence of a positive spectral gap $|\lambda_1|$
of $\mathcal{L}.$
\begin{figure}
   \epsfysize=2.5cm
 $$\epsfbox{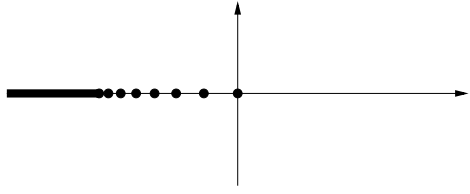}$$
 \caption{Spectrum of the collision operator $\Q$ in $\mathcal{H}$.}
\end{figure}
Note that such a result, combined with Proposition \ref{negat}, has
important consequence on the entropy production, since it can be
shown in an easy way that the $H$-Theorem reads as
$$\dfrac{\d }{\d t}\left\|f(t)-\mathcal{M}\right\|_{\mathcal{H}}^2
=\langle \mathcal{L}f(t),f(t)\rangle_{\mathcal{H}}.$$ Consequently,
the {\it Dirichlet form} $B(f)=\langle \mathcal{L}f ,f
\rangle_{\mathcal{H}}$ plays the role of entropy-dissipation
functional and the existence of a spectral gap $|\lambda_1|$ is
equivalent to the following coercivity estimate for $B(f)$:
$$B(f) \geq -|\lambda_1|\, \|f\|_{\mathcal{H}}^2 \qquad \qquad \forall f \perp \mathrm{span}(\mathcal{M}).$$
One deduces easily the following corollary on the exponential trend towards equilibrium:
\begin{cor} Let $f_0(\v) \in \mathcal{H}$ and let $f(t,\v)$ be the unique solution to the linear homogeneous Boltzmann equation
\eqref{homogeneous}. Then, there is some constant $C \geq 0$ such that
$$\left\|f(t,\cdot)-\mathcal{M}\right\|_{\mathcal{H}} \leq C\, \exp(-|\lambda_1| t)\left\|f_0-\mathcal{M}\right\|_{\mathcal{H}},\qquad \text{ for any } \quad t \geq 0,$$
where $0 < |\lambda_1| \leq \nu_0$ is provided by Theorem \ref{spectral}.
\end{cor}

 We refer the reader to \cite{moto} for details on the
matter, and in particular, for an {\it explicit estimate of the spectral gap} $|\lambda_1|$.

\section{Application to the honest solutions  of the Boltzmann equation}\label{sec:honest}
\subsection{Conservative solutions} We are interested in this section in applying the result of the
previous section to prove the existence of honest solutions to the
linear Boltzmann equation for dissipative hard--spheres
\begin{equation}\label{bolt} {\partial_t}f(\x,\v,t)+\v \cdot \nabla_\x
f(\x,\v,t)=\Q(f)(\x,\v,t),\end{equation} with initial condition
$$f(\x,\v,0)=f_0(\x,\v) \in L^1(\R \times \R,\d\x \otimes \d\v),$$ where the collision operator $\Q$ is given by Eq. \eqref{linearoperator}. Recall that
the streaming operator $T_0$ is defined by
$$\D(T_0)=\{f \in \mathcal{X}\,,\v \cdot \nabla_\x f \in \mathcal{X}\}, \qquad T_0f=-\v
\cdot \nabla_\x f$$ where $\mathcal{X}=L^1(\R \times \R,\d\x\otimes\d\v)$.  One can
define then the multiplication operator $\Sigma$ by
$$\D(\Sigma)=\{f \in \mathcal{X}\,,\sigma\,f \in \mathcal{X}\},\qquad \Sigma f(x,v)=-\sigma(v)\,f(x,v)$$
where, as  in the previous Section, $\sigma(\v)$ is the collision
frequency corresponding to {\it dissipative hard spheres}
interactions and given by Eq. \eqref{sigmaHS}.  The following generation
result is well-known \cite{arloban}
\begin{theo}

The operator $T_0$ generates a $C_0$-semigroup of isometries $\uot$
of $\mathcal{X}$ given by:
$$U(t)f(\x,\v)=f(\x-t\v,\v), \qquad t \geq 0.$$
The operator $A=T_0+\Sigma$ with domain $\D(A)=\D(T_0) \cap
\D(\Sigma)$ is the generator of a contractions $C_0$-semigroup $\vt$
given by
$$V(t)f(\x,\v)=\exp(-\sigma(\v) t)f(\x-t\v,\v),\qquad t \geq 0.$$
\end{theo}

Let us define now $K$ as the gain operator $\Q^+$ endowed with the
domain of $A$:
$$\D(K)=\D(A),\qquad Kf(\x,\v)=\Q^+(f)(\x,\v)=\epsilon^{-2}\It |\q \cdot \n|f(\x,\v_{\star})\mathcal{M}_1(\w_{\star})\d \w \d
\n.$$
 It is clear from \eqref{conserva} that, for any $f \in \D(K)$,
\begin{equation}\label{conserv}\int_{\R \times \R}(Af + Kf)\d \x\d \v=0,\end{equation} while $Kf
\geq 0$ for any $f \in \D(K)$, $f \geq 0.$ Then, the following
generation result is a direct consequence of \cite{Ar,voigt}:
\begin{theo} There exists a positive contractions semigroup $\zt$ in $\mathcal{X}$ whose
generator $G$ is an extension of $A+K$. Moreover, $\zt$ is minimal,
i.e. if $(T(t))_{t \geq 0}$ is a positive $C_0$-semigroup generated
by an extension of $A+K$, then $T(t) \geq Z(t)$ for any $t \geq
0.$\end{theo}

The natural question is now to determine whether the "formal" mass
conservation identity \eqref{massconser} can be made rigorous.
Namely, one aims to prove that, for any nonnegative $f \in
\mathcal{X}$, the following holds:
$$\|Z(t)f\|=\|f\|,\qquad \forall t >0.$$
The important point to be noticed is the following. If
$G=\overline{A+K}$, then any function $\varphi \in \D(G)$ can be
approximated by a sequence $(\varphi_n)_n \subset \D(A+K)=\D(A)$
such that $\varphi_n \to \varphi$ and $(A+K)\varphi_n \to G\varphi$
as $n \to \infty.$ In particular, \eqref{conserv} implies that
$$\int_{\R  \times \R}G\varphi\,\d \x \d \v=\lim_{n \to \infty}\int_{\R \times
\R}(A+K)\varphi_n \d \x \d \v=0, \qquad \forall \varphi \in \D(G).$$
Now, for any given initial datum $f_0 \in \D(G)$, $f_0 \geq 0$, the
solution $f(t)=Z(t)f_0$ of \eqref{bolt} is such that
$$\dfrac{\d}{\d t}\|f(t)\|=\int_{\R \times \R}\dfrac{\d }{\d t}f(t) \d \x \d \v=\int_{\R \times \R}Gf(t) \d \x \d \v=0,$$
i.e.
$$\|f(t)\|=\|f_0\|, \qquad \forall t \geq 0.$$
 This means that, if $G=\overline{A+K}$, then the solutions to
the linear Boltzmann equation \eqref{bolt} are conservative. On the other hand,
if $G$ is a larger extension of $A+K$ than
$\overline{A+K}$, then there may be a loss of particles in the
evolution (see \cite{arloban}
for the matter as well as \cite{labl} for examples of transport equation for which such a loss of particles occurs because of boundary conditions). Precisely, if $G \neq \overline{A+K}$ then there exists
$f_0 \in \mathcal{X}$, $f_0 \geq 0$ such that
$$\|Z(t)f_0\| < \|f_0\| \qquad \text{ for some } t > 0.$$

This shows that the determination of the domain $\D(G)$ of $G$ is of
primary importance in the study of the Boltzmann equation. This is
the main concern of the so--called \textit{substochastic
perturbation theory of $C_0$-semigroups}  \cite{arloban}.

We point out that  the question   of the honesty of the semigroup
governing the Boltzmann equation also arises in the study of the
{\it space-homogeneous} version of the latter equation. Indeed, it
is the unboundedness of the collision frequency (and consequently
that of whole collision operator $\Q$) that may give rise to {\it
dishonest solutions} to the Boltzmann equation. Actually, to prove
the honesty of the $C_0$-semigroup $\zt$, we will adopt the strategy
developed first in \cite{Ar} and systematized in \cite{arloban}.
More precisely, we will show that the gain operator $K$ fullfils the
assumption of \cite{Ar}:
\begin{propo} There exists $C>0$ such that, for any fixed
$\varrho >0$,
$$\esup_{|\v-\u| \leq \varrho} \int_{|\v'-\u| \geq \varrho}k(\v',\v)\d\v' \leq C.$$
\end{propo}
\begin{proof} Since the
kernel $k(\v,\v')$ differs from the corresponding one for classical
gas, except from numerical constants, one can apply {\it mutatis
mutandis} the technical calculations of \cite[Section 4.1]{Ar}
(see also \cite[p. 329-330]{arloban}) to get the desired estimate.\end{proof}

As a consequence, one deduces immediately from \cite{Ar}, the main
result of this section:
\begin{theo}\label{main}
The generator $G$ of the minimal semigroup $\zt$ is given by
$$G=\overline{A+K}.$$ In particular, the $C_0$-semigroup $\zt$ is
honest and
$$\int_{\R \times \R}Z(t)f(\x,\v)\d\x\d\v=\int_{\R \times \R}f(\x,\v)\d\x\d\v \qquad \text{ for any } f \in \mathcal{X} \: \text{
and any } t \geq 0.$$\end{theo}

\subsection{Consequence on the entropy production}

As a direct application of the above result (Theorem \ref{main}), we
give a direct rigorous  proof of the linear $H$--Theorem of
\cite{loto}. In order to stay in the formalism of \cite{loto}, we
shall restrict ourselves to the space-homogenous case. Precisely,
let $Y$ denote the set of functions depending only on the velocity
and integrable with respect to velocities:
$$Y=L^1(\R,\d\v),$$
equipped with its natural norm $\|\cdot\|_Y.$ For any nonnegative
$f$ and $g$ in $Y$, we define the  \textit{information}  of $f$
with respect to $g$ by
$$\mathcal{I}(f|g)=\int_{ \R} \bigg(f(\v) \ln f(\v) -f(\v)\ln g(\v)\bigg) \d \v$$
with the conventions $0\ln 0=0$ and $x\ln 0=-\infty$ for any $x
>0$. This means that the information is nothing but the entropy functional $H_{\Phi}$ for the particular choice of $\Phi(s)=s \ln s.$
One recalls the main result of  \cite{voigtentropy}:
\begin{theo}\label{voigt}
Let $U$ be a stochastic operator of $\,Y$, i.e. $U$ is a positive operator such that $\|Uf\|_Y=\|f\|_Y$ for any $f \in Y,$ $f
\geq 0$. Then,
$$\mathcal{I}(Uf|Ug) \leq \mathcal{I}(f|g)$$
for any nonnegative $f$, $g$ in $Y$. In particular, if $g \in Y$ is
a nonnegative fixed  point of $U$ then,
$$\mathcal{I}(Uf|g) \leq \mathcal{I}(f|g), \qquad \forall f \in Y,\:f \geq 0.$$
\end{theo}
According to the results of the previous section, it is not
difficult to see that the restriction of $\zt$ to $Y$ is a
$C_0$-semigroup of stochastic operators of $Y$.  Since the unique
equilibrium state $\mathcal{M} \in Y$ is space independent, one sees
that $(T_1+K)\mathcal{M}=0$ and, in particular,
$$Z_Y(t)\mathcal{M}=\mathcal{M}, \qquad \qquad \forall t \geq 0.$$
Combining this with Theorem \ref{voigt}, one obtains a rigorous and
direct proof of the $H$--Theorem \cite[Theorem 5.1]{loto}:
\begin{theo}
Let $f_0 \in Y$ be a given nonnegative (space homogeneous)
distribution function with unit mass, i.e. $\|f_0\|_Y=1$. Assume
that $\mathcal{I}(f_0|\mathcal{M}) < \infty,$ then
$$\dfrac{\d}{\d t}\mathcal{I}(f(t)|\mathcal{M}) \leq 0, \qquad \qquad (t \geq 0),$$
where $f(t)=Z_Y(t)f_0=Z(t)f_0$ is the unique solution to
\eqref{bolt} with $f(0)=f_0$.
\end{theo}
\begin{nb} Note that, once the conservativity of the solution to the
Boltzmann equation asserted by Theorem \ref{main}, the above
$H$-Theorem can be proved by usual standard method of kinetic
theory. However, we insist on the fact that such standard proofs
require the solution $f(t,\v)$ to be conservative and, in some
sense, the use of the substochastic semigroup theory.\end{nb}


{\small
\begin{thebibliography}{20}
\bibitem{Ar}
{\sc L. Arlotti,} A perturbation theorem for positive contraction
semigroups on $L\sp 1$-spaces with applications to transport
equations and Kolmogorov's differential equations. {\em Acta Appl.
Math.}  {\bf 23},  129--144, 1991.

\bibitem{labl}
\textsc{L. Arlotti, B. Lods}, Substochastic semigroups for transport
equations with conservative boundary conditions, \textit{J. Evol.
Equations} {\bf 5}  485--508,  2005.

\bibitem{arloban}
\textsc{J. Banasiak, L. Arlotti}, \textbf{Perturbations of positive
semigroups with applications}, Springer-Verlag, 2005.


\bibitem{poschel}
\textsc{N. V. Brilliantov, T. P\"{o}schel}, \textbf{Kinetic theory of
granular gases}, Oxford University Press, 2004.

\bibitem{pareschi}
\textsc{S.Brull, L.Pareschi,} Dissipative hydrodynamic models for
the diffusion of impurities in a gas,  \textit{Appl. Math. Lett.}
\textbf{19} 516--521, 2006.

\bibitem{carleman}
{\sc T. Carleman,}
\newblock Probl\`emes math\'ematiques dans la th\'eorie cin\'etique des gaz.
\newblock {\it Publications Scientifiques de l'Institut
Mittag-Leffler,} {\bf 2}, 1957.
\bibitem{cerci}
\textsc{C. Cercignani, R. Illner and  M. Pulvirenti}, \textbf{The
mathematical theory of dilute gases}, Springer--Verlag, New York,
1994.
\bibitem{ernst}
\textsc{M. H. Ernst, R. Brito.}
\newblock{Scaling solutions of inelastic Boltzmann equation with over--populated high energy
tails.}
\newblock {\em J. Statist. Phys.} {\bf 109}, 407--432, 2002.

\bibitem{ferrari}
\textsc{E. Ferrari, L. Pareschi}, Modelling and numerical methods
for the diffusion of impurities in a gas, \textit{Int. J. Numer.
Meth. Fluids}, to appear.

\bibitem{garzosurvey}
\textsc{V. Garz\`{o},} Kinetic Theory for Binary Granular Mixtures at
Low-Density, preprint, \texttt{ArXiv:0704.1211}.

\bibitem{garzo}
\textsc{V. Garz\`{o}, J. M. Montanero,} Diffusion of impurities in a
granular gas, \textit{Physical Review E}, \textbf{69},  2004.

\bibitem{Grad}
{\sc H. Grad,}
\newblock Asymptotic theory of the {B}oltzmann equation. {II}.
\newblock {\it Rarefied Gas Dynamics (Proc. 3rd Internat. Sympos., Palais de
l'UNESCO, Paris, 1962)}, Vol. I, 26--59, 1963.
\bibitem{kuscer}
{\sc I. Kuscer, M. M. R. Williams}, Relaxation constants of a
uniform hard-sphere gas, {\it  Physics of Fluids}, {\bf 10},
1922--1927, 1967.
\bibitem{loto} {\sc B. Lods, G. Toscani,}
\newblock{The dissipative linear Boltzmann equation for hard spheres,}
\newblock {\em J. Statist. Phys.} \textbf{117}, 635--664, 2004.
\bibitem{Mo}
{\sc F. A. Molinet,} Existence, uniqueness and properties of the
solutions of the Boltzmann kinetic equation for a weakly ionized
gas, {\em J. Math. Phys.} {\bf 18},  984--996, 1977.
\bibitem{piasecki}
\textsc{P. A. Martin, J. Piaceski,} Thermalization of a particle by
dissipative collisions, \emph{Europhys. Lett.} \textbf{46},
 613--616, 1999.

\bibitem{mischler}
{\sc S. Mischler, C. Mouhot,} Cooling process for inelastic
Boltzmann equations for hard spheres, Part II: Self-similar
solutions and tail behavior, {\em J. Statist. Phys.} \textbf{124},
703--746, 2006.

\bibitem{strain}
{\sc C. Mouhot, R. Strain,}  Spectral gap and coercivity estimates
for the linearized Boltzmann collision operator without angular
cutoff,  {\em J. Maths Pures Appl.},  in press.


\bibitem{moto}
{\sc C. Mouhot, G. Toscani}, Relaxation rate and diffusive limit for
inelastic scattering Boltzmann models, Work in progress.
\bibitem{rolf}
{\sc R. Pettersson,}
\newblock{On solutions to the linear Boltzmann
equation for granular gases,} \newblock{\em Transp. Theory Statist.
Phys.} \textbf{33},  527--543, 2004.
\bibitem{poupaud}
{\sc F. Poupaud,} Runaway phenomena and fluid approximation under
high fields in semiconductor kinetic theory, {\em Z. Angew. Math.
Mech.} \textbf{72},  359--372, 1992.

\bibitem{spiga}
{\sc G. Spiga, G. Toscani,} The dissipative linear Boltzmann
equation, {\em Appl. Math. Lett.} {\bf 17}, 255--301, 2004.
\bibitem{villani}
{\sc C. Villani,}
\newblock{Mathematics of granular materials,}
\newblock{\em J. Statist. Phys.} {\bf 124},  781--822, 2006.

\bibitem{voigtentropy}
\textsc{J. Voigt}, Stochastic operators, information, and entropy,
\textit{Comm. Math. Phys.} \textbf{81},  31--38, 1981.

\bibitem{voigt}
\textsc{J. Voigt}, On substochastic $C_0$-semigroups and their
generators, \textit{Transp. Theory Stat. Phys.} {\bf 16}, 453--466,
1987.


\end{thebibliography}
}
 \end{document}